\documentclass{amsart}
\usepackage{enumerate}
\usepackage{amssymb}
\usepackage{graphicx}
\usepackage{url}

\author{Bjarke Hammersholt Roune}

\title{The Slice Algorithm For Irreducible Decomposition of Monomial Ideals}

\theoremstyle{definition}
\newtheorem{theorem}{Theorem}
\newtheorem{lemma}[theorem]{Lemma}
\newtheorem{definition}[theorem]{Definition}
\newtheorem{proposition}[theorem]{Proposition}
\newtheorem{corollary}[theorem]{Corollary}

\newtheorem{example}[theorem]{Example}

\DeclareMathOperator{\lcm}{lcm}
\DeclareMathOperator{\irr}{irr}
\DeclareMathOperator{\supp}{supp}
\DeclareMathOperator{\msmSym}{msm}
\DeclareMathOperator{\conSym}{con}

\def\cocoa{{\hbox{\rm C\kern-.13em o\kern-.07em C\kern-.13em o\kern-.15em A}}}

\newcommand{\R}{\mathbb{R}}
\newcommand{\Q}{\mathbb{Q}}
\newcommand{\Z}{\mathbb{Z}}
\newcommand{\N}{\mathbb{N}}

\newcommand{\p}{\ensuremath{^\prime}}

\newcommand{\biimp}{\Leftrightarrow}
\newcommand{\imp}{\Rightarrow}

\newcommand{\setBuilder}[2]{\left\{#1\left|#2\right.\right\}}
\newcommand{\idealBuilder}[2]{\left\langle#1\left|#2\right.\right\rangle}

\newcommand{\card}[1]{\left|#1\right|}

\newcommand{\proofPart}[1]{{\bf $\boldsymbol{#1}$:}}

\newcommand{\msm}[1]{\msmSym\left({#1}\right)}
\newcommand{\con}[3]{\conSym\left({#1},{#2},{#3}\right)}
\newcommand{\cont}[1]{\conSym\left({#1}\right)}
\newcommand{\decom}[1]{\irr\left({#1}\right)}
\newcommand{\ming}[1]{\min\left({#1}\right)}

\newcommand{\ideal}[1]{\left<#1\right>}

\newcommand{\set}[1]{\left\{{#1}\right\}}

\newcommand{\degx}[2]{\ensuremath{\deg_{x_{#1}}\!\!\left({#2}\right)}}

\newcommand{\projSym}{\ensuremath{\pi}}
\newcommand{\proj}[1]{\ensuremath{\projSym\left({#1}\right)}}

\begin{document}

\begin{abstract}
Irreducible decomposition of monomial ideals has an increasing number
of applications from biology to pure math. This paper presents the
Slice Algorithm for computing irreducible decompositions, Alexander
duals and socles of monomial ideals. The paper includes experiments
showing good performance in practice.
\end{abstract}
\today
\maketitle

\setcounter{tocdepth}{1}
\tableofcontents 

\section{Introduction}

The main contribution of this paper is the Slice Algorithm, which is
an algorithm for the computation of the irreducible decomposition of
monomial ideals. To irreducibly decompose an ideal is to write it as
an irredundant intersection of irreducible ideals.

Irreducible decomposition of monomial ideals has an increasing number
of applications from biology to pure math. Some examples of this are
the Frobenius problem \cite{bfrob, frobPoint}, the integer programming
gap \cite{igap}, the reverse engineering of biochemical networks
\cite{bionet}, tropical convex hulls \cite{tropicalHull}, tropical
cyclic polytopes \cite{tropicalHull}, secants of monomial ideals
\cite{secant}, differential powers of monomial ideals \cite{diffPower}
and joins of monomial ideals \cite{secant}.

Irreducible decomposition of a monomial ideal $I$ has two
computationally equivalent guises. The first is as the \emph{Alexander
dual} of $I$ \cite{alexdual}, and indeed some of the references above
are written exclusively in terms of Alexander duality rather than
irreducible decomposition. The second is as the \emph{socle} of the
vector space $R/I\p$, where $R$ is the polynomial ring that $I$
belongs to and $I\p:=I+\ideal{x_1^t,\ldots,x_n^t}$ for some integer
$t>>0$. The socle is central to this paper, since what the Slice
algorithm actually does is to compute a basis of the socle.

Section \ref{sec:prelim} introduces some basic notions we will need
throughout the paper and Section \ref{sec:basicAlg} describes an
as-simple-as-possible version of the Slice Algorithm. Section
\ref{sec:improve} contains improvements to this basic version of the
algorithm and Section \ref{sec:selection} discusses some heuristics
that are inherent to the algorithm. Section \ref{sec:apps} examines
applications of irreducible decomposition, and it describes how the
Slice Algorithm can use bounds to solve some optimization problems
involving irreducible decomposition in less time than would be needed
to actually compute the decomposition. Finally, Section
\ref{sec:bench} explores the practical aspects of the Slice Algorithm
including benchmarks comparing it to other programs for irreducible
decomposition.

The Slice Algorithm was in part inspired by an algorithm for
Hilbert-Poincar\'e series due to Bigatti Et
Al. \cite{bigattiEtAlHilbSerlg}. The Slice Algorithm generalizes
versions of the staircase-based algorithm due to Gao and Zhu
\cite{artinianStairIrr} (see Section \ref{sssec:purePowerSelect}) and
the Label Algorithm due to Roune \cite{rouneLabelAlg} (see Section
\ref{sssec:labelSelect}).

\section{Preliminaries}
\label{sec:prelim}

This section briefly covers some notation and background on monomial
ideals that is necessary to read the paper. We assume throughout the
paper that $I$, $J$ and $S$ are monomial ideals in a polynomial ring
$R$ over some arbitrary field $\kappa$ and with variables
$x_1,\ldots,x_n$ where $n\geq 2$. We also assume that $a$, $b$, $p$,
$q$ and $m$ are monomials in $R$. When presenting examples we use the
variables $x$, $y$ and $z$ in place of $x_1$, $x_2$ and $x_3$ for
increased readability.

\subsection{Basic Notions From Monomial Ideals}

If $v\in\N^n$ then $x^v:=x_1^{v_1}\cdots x_n^{v_n}$. We define $\sqrt
{x^v}:=x^{\supp(v)}$ where $(\supp(v))_i:=\min(1,v_i)$. Define $\proj
m:=\frac{m}{\sqrt{m}}$ such that
e.g. $\proj{x^{(0,1,2,3)}}=x^{(0,0,1,2)}$.

The rest of this section is completely standard. A \emph{monomial
ideal} $I$ is an ideal generated by monomials, and $\ming I$ is the
unique minimal set of monomial generators. The ideal $\ideal M$ is the
ideal generated by the elements of the set $M$. The \emph{colon ideal}
$I:p$ is defined as $\ideal{m|mp\in I}$.

An ideal $I$ is \emph{artinian} if there exists a $t\in\N$ such that
$x_i^t\in I$ for $i=1,\ldots,n$. A monomial of the form $x_i^t$ is a
\emph{pure power}. A monomial ideal is \emph{irreducible} if it is
generated by pure powers. Thus $\ideal{x^2,y}$ is irreducible while
$\ideal{x^2y}$ is not. Note that $\ideal{x}\subseteq\kappa[x,y]$ is
irreducible and not artinian.

Every monomial ideal $I$ can be written as an irredundant intersection
of irreducible monomial ideals, and the set of ideals that appear in
this intersection is uniquely given by $I$. This set is called the
\emph{irreducible decomposition} of $I$, and we denote it by $\decom
I$. Thus $\decom{\ideal{x^2,xy,y^3}}=
\set{\ideal{x^2,y},\ideal{x,y^3}}$.

The \emph{radical} of a monomial ideal $I$ is
$\sqrt{I}:=\idealBuilder{\sqrt{m}}{m\in\ming I}$.  A monomial ideal
$I$ is \emph{square free} if $\sqrt I=I$. A monomial ideal is
(strongly) \emph{generic} if no two distinct elements of $\ming I$
raise the same variable $x_i$ to the same non-zero power
\cite{scarfComplex, genericMonomial}. Thus $\ideal{x^2y,xy^2}$ is
generic while $\ideal{xyz^2,xy^2z}$ is not as both minimal generators
raise $x$ to the same power. In this paper we informally talk of a
monomial ideal being more or less generic according to how many
identical non-zero exponents there are in $\ming I$.

\begin{figure}[!ht]
\centering
\begin{tabular}{c|c|c}
\includegraphics[width=0.3\textwidth]{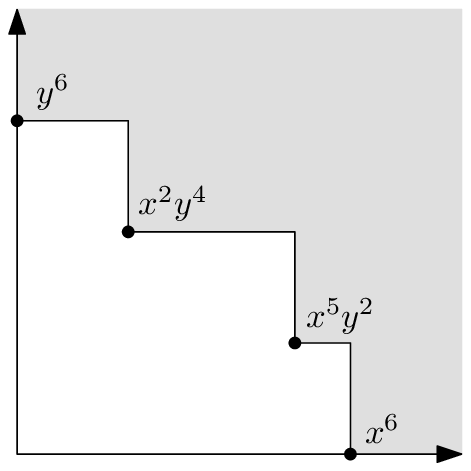}&
\includegraphics[width=0.3\textwidth]{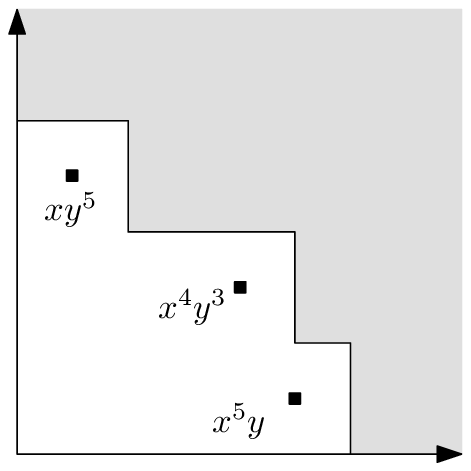}&
\includegraphics[width=0.3\textwidth]{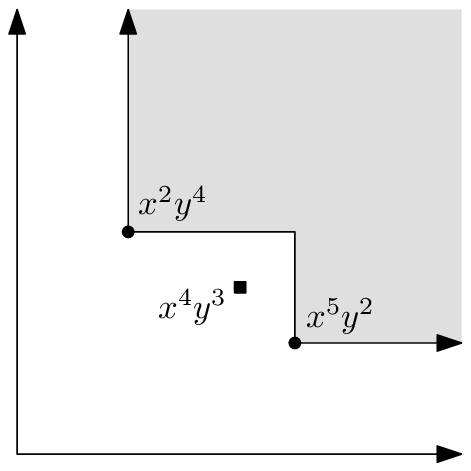}\\
(a)&(b)&(c)\\
\end{tabular}
\caption{Examples of monomial ideals.}
\label{fig:msm}
\end{figure}

A \emph{standard monomial} of $I$ is a monomial that does not lie
within $I$. The \emph{exponent vector} $v\in\N^n$ of a monomial $m$ is
defined by $m=x^v=x_1^{v_1}\cdots x_n^{v_n}$. Define $\degx i
{x^v}:=v_i$. We draw pictures of monomial ideals in 2 and 3
dimensions by indicating monomials by their exponent vector and
drawing line segments separating the standard monomials from the
non-standard monomials. Thus Figure \ref{fig:msm}(a) displays a
picture of the monomial ideal $\ideal{x^6,x^5y^2,x^2y^4,y^6}$.

\subsection{Maximal Standard Monomials, Socles And Decompositions}
\label{ssec:msm}

In this section we look into socles and their relationship with
irreducible decomposition. We also note the well known fact that the
maximal standard monomials of $I$ form a basis of the socle of $R/I$.

Given the generators $\ming I$ of a monomial ideal $I$, the Slice
Algorithm computes the maximal standard monomials of $I$. We will need
some notation for this.

\begin{definition}[Maximal standard monomial]
A monomial $m$ is a \emph{maximal standard monomial} of $I$ if
$m\notin I$ and $mx_i\in I$ for $i=1,\ldots,n$. The set of maximal
standard monomials of $I$ is denoted by $\msm I$.
\end{definition}

The \emph{socle} of $R/I$ is the vector space of those $m\in R/I$ such
that $mx_i=0$ for $i=1,\ldots,n$. It is immediate that
$\setBuilder{m+I}{m\in\msm I}$ is a basis of this socle.

\begin{example}
Let $I:=\ideal{x^6,x^5y^2,x^2y^4,y^6}$ be the ideal in Figure
\ref{fig:msm}(a). Then $\msm I=\set{x^5y,x^4y^3,xy^5}$ as indicated on
Figure \ref{fig:msm}(b). Let $J:=\ideal{x^5y^2,x^2y^4}$.  Then
$\msm{J}=\set{x^4y^3}$ as indicated in Figure
$\ref{fig:msm}(c)$. Finally, $\msm{\ideal{x^5y^2}}=\emptyset$.
\end{example}

We will briefly describe the standard technique for obtaining $\decom
I$ from $\msm I$ \cite{scarfComplex}. Choose some integer $t>>0$ and
define $\phi(x^m)=\idealBuilder{x_i^{m_i+1}}{m_i+1<t}$.

\begin{proposition}[\protect{\cite[ex. 5.8]{cca}}]
The map $\phi$ is a bijection from
$\msm{I+\ideal{x_1^t,\ldots,x_n^t}}$ to $\decom I$.
\end{proposition}

\begin{example}
Let $I:=\ideal{x^2,xy}$ and $I\p:=I+\ideal{x^t,y^t}=\ideal{x^2,xy,y^3}$
where $t=3$. Then $\msm{I\p}=\set{x, y^2}$ which $\phi$ maps to
$\set{\ideal{x^2,y},\ideal{x}}=\decom I$.
\end{example}

\subsection{Labels}

We will have frequent use for the notion of a label.

\begin{definition}[$x_i$-label]
Let $d$ be a standard monomial of $I$ and let $m\in\ming I$. Then $m$
is an \emph{$x_i$-label} of $d$ if $m|dx_i$.
\end{definition}

Note that if $m$ is an $x_i$-label of $d$, then $\degx i m=\degx i
d+1$. Also, a standard monomial $d$ is maximal if and only if it has
an $x_i$-label $m_i$ for $i=1,\ldots,n$. So in that case $dx_1\cdots
x_n=\lcm_{i=1}^n m_i$.

\begin{figure}[!ht]
\centering
\begin{tabular}{c|c}
\includegraphics[width=0.3\textwidth]{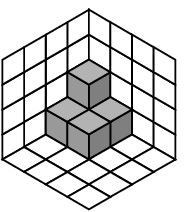}&
\includegraphics[width=0.3\textwidth]{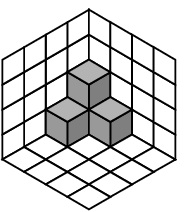}\\
(a)&(b)\\
\end{tabular}
\caption{Examples of monomial ideals.}
\label{fig:label}
\end{figure}

\begin{example}
Let $I:=\ideal{x^2,xz,y^2,yz,z^2}$ be the ideal in Figure
\ref{fig:label}(a). Then the maximal standard monomials of $I$ are
$\msm I=\set{xy,z}$. We see that $z$ has $xz$ as an $x$-label, $yz$ as
a $y$-label and $z^2$ as a $z$-label. Also, $xy$ has $x^2$ as an
$x$-label and $y^2$ as a $y$-label, while it has both of $xz$ and $yz$
as $z$-labels.

Let $J:=I+\ideal{xy}$ be the ideal in Figure \ref{fig:label}(b). Then
$\msm I=\set{x,y,z}$. Note that even though $xy$ divides $z\cdot xyz$,
it is not a label of $z$, because it does not divide $z\cdot x$,
$z\cdot y$ or $z\cdot z$.
\end{example}

\section{The Slice Algorithm}
\label{sec:basicAlg}

In this section we describe a basic version of the Slice
Algorithm. The Slice Algorithm computes the maximal standard monomials
of a monomial ideal given the minimal generators of that ideal.

A fundamental idea behind the Slice Algorithm is to consider certain
subsets of $\msm I$ that are represented as \emph{slices}. We will
define the meaning of the term slice shortly. The algorithm starts out
by considering a slice that represents all of $\msm I$. It then
processes this slice by splitting it into two simpler slices. This
process continues recursively until the slices are simple enough that
it is easy to find any maximal standard monomials within them.

From this description, there are a number of details that need to be
explained. Section \ref{ssec:slices} covers what slices are and how to
split them while Section \ref{ssec:base} covers the base case. Section
\ref{ssec:term} proves that the algorithm terminates and Section
\ref{ssec:pseudo} contains a simple pseudo-code implementation of the
algorithm.

\subsection{Slices And Splitting}
\label{ssec:slices}
In this section we explain what slices are and how to split them. We
start off with the formal definition of a slice and its content.

\begin{definition}[Slice and content]
A \emph{slice} is a 3-tuple $(I,S,q)$ where $I$ and $S$ are monomial
ideals and $q$ is a monomial. The \emph{content} of a slice is defined
by $\con I S q:=(\msm I\setminus S)q$.
\end{definition}

Example \ref{ex:split} shows how this definition is used.

\begin{figure}[!ht]
\centering
\begin{tabular}{c|c|c}
\includegraphics[width=0.3\textwidth]{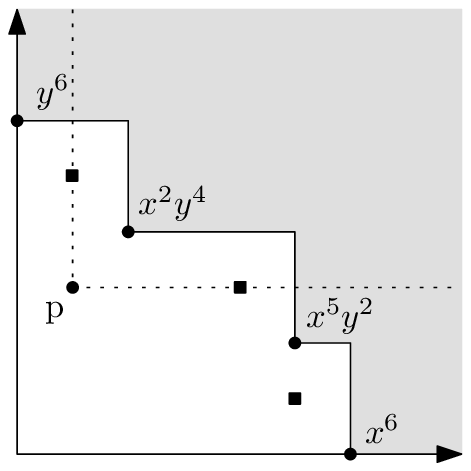}&
\includegraphics[width=0.3\textwidth]{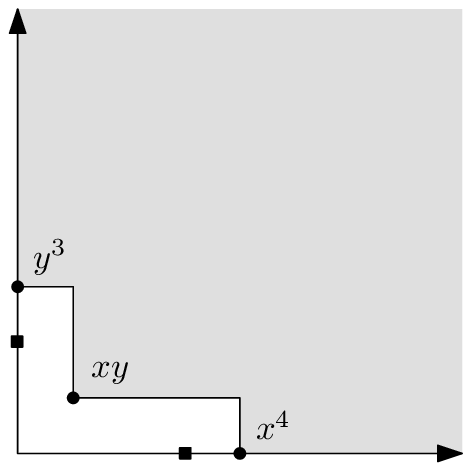}&
\includegraphics[width=0.3\textwidth]{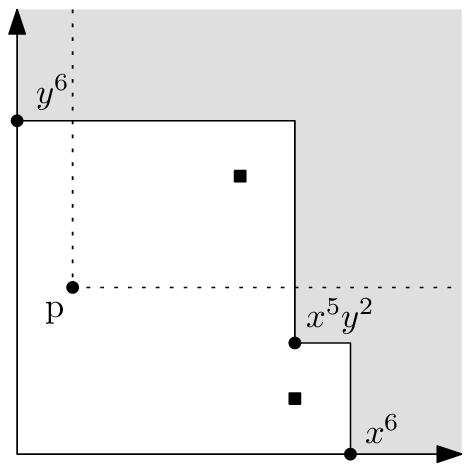}\\
(a)&(b)&(c)\\
\end{tabular}
\caption{Illustrations for example \ref{ex:split}.}
\label{fig:split}
\end{figure}

\begin{example}
\label{ex:split}
Let $I:=\ideal{x^6,x^5y^2,x^2y^4,y^6}$ and $p:=xy^3$. Then $I$ is the
ideal depicted in Figure \ref{fig:split}(a), where $\ideal p$ is
indicated by the dotted line and $\msm I=\set{x^5y,x^4y^3,xy^5}$ is
indicated by the squares. We will compute $\msm I$ by performing a
step of the Slice Algorithm.

Let $I_1$ be the ideal $I:p=\ideal{y^3,xy,x^4}$, as depicted in Figure
\ref{fig:split}(b), where $\msm{I_1}=\set{x^3,y^2}$ is indicated by
the squares. As can be seen by comparing figures \ref{fig:split}(a)
and \ref{fig:split}(b), the ideal $I_1$ corresponds to the part of the
ideal $I$ that lies within $\ideal p$. Thus it is reasonable to expect
that $\msm{I_1}$ corresponds to the subset of $\msm I$ that lies
within $\ideal p$, which turns out to be true, since
\begin{equation}
\label{eqn:splitEx1}
\msm{I_1}p=\set{x^4y^3,xy^5}=\msm I\cap\ideal p
\end{equation}
It now only remains to compute $\msm{I}\setminus\ideal p$. Let
$I_2:=\ideal{x^6,x^5y^2,y^6}$ as depicted in Figure
\ref{fig:split}(c), where $\msm{I_2}:=\set{x^5y,x^4y^5}$ is indicated
by the squares. The dotted line indicates that we are ignoring
everything inside $\ideal p$. It happens to be that one of the minimal
generators of $I$, namely $x^2y^4$, lies in the interior of $\ideal
p$, which allows us to ignore that minimal generator. We are looking
at $I_2$ because
\begin{equation}
\label{eqn:splitEx2}
\msm{I_2}\setminus\ideal p=\set{x^5y}=\msm I\setminus\ideal p
\end{equation}
By combining Equation \eqref{eqn:splitEx1} and Equation
\eqref{eqn:splitEx2}, we can compute $\msm I$ in terms of $\msm{I_1}$,
$\msm{I_2}$ and $p$.

Using the language of slices, we have split the slice $A:=(I,\ideal
0,1)$ into the two slices $A_1:=(I_1,\ideal 0,p)$ and
$A_2:=(I_2,\ideal p,1)$. By Equations \eqref{eqn:splitEx1} and
\eqref{eqn:splitEx2}, we see that $\cont{A_1}=\msm I\cap\ideal p$ and
$\cont{A_2}=\msm I\setminus\ideal p$. Thus
\[
\cont A=\msm I=\cont{A_1}\cup\cont{A_2}
\]
where the union is disjoint.
\end{example}

Having defined slices and their content, we can now explain how to
\emph{split} a slice into two smaller slices. This is done by choosing
some monomial $p$, called the \emph{pivot}, and then to consider the
following trivial equation.
\begin{equation}
\label{eqn:splitTrivial}
\con I S q = \Big(\con I S q \cap\ideal{qp}\Big)
\cup\Big( \con I S q\setminus\ideal{qp}\Big)
\end{equation}
The idea is to express both parts of this disjoint union as the
content of a slice. This is easy to do for the last part, since
\[
\label{eqn:sliceIbyB}
\con I S q\setminus\ideal{qp}=\con I {S+\ideal p} q
\]
Expressing the first part of the union as the content of a slice can
be done using the following equation, which we will prove at the end
of this section.
\[
\msm I\cap\ideal p=\msm{I:p}p
\]
which implies that (see Example \ref{ex:split})
\[
\label{eqn:sliceIbyM}
\con I S q\cap\ideal{qp}=\con {I:p} {S:p} {qp}
\]
Thus we can turn Equation \eqref{eqn:splitTrivial} into the following.
\begin{equation}
\label{eqn:split}
\con I S q=\con {I:p} {S:p} {qp}\,\,\cup\,\, \con I {S+\ideal p} q
\end{equation}

Equation \eqref{eqn:split} is the basic engine of the Slice
Algorithm. We will refer to it and its parts throughout the paper, and
we need some terminology to facilitate this. The process of applying
Equation \eqref{eqn:split} is called a \emph{pivot split}. We will
abbreviate this to just \emph{split} when doing so should not cause
confusion.

Equation \eqref{eqn:split} mentions three slices, and we give each of
them a name. We call the left hand slice $(I,S,q)$ the \emph{current
slice}, since it is the slice we are currently splitting. We call the
first right hand slice $(I:p,S:p,qp)$ the \emph{inner slice}, since
its content is inside $\ideal{qp}$, and we call the second right hand
slice $(I,S+\ideal p,q)$ the \emph{outer slice}, since its content is
outside $\ideal{qp}$.

It is not immediately obvious why it is easier to compute the outer
slice's content $\con I {S+\ideal p} q$ than it is to compute the
current slice's content $\con I S q$. The following equation shows how
it can be easier. See Proposition \ref{prop:split} for a proof.
\begin{equation}
\label{eqn:normalize}
\msm I\setminus S=\msm{I\p}\setminus S,
\,\,I\p:=\idealBuilder{m\in\ming I}{\proj m\notin S}
\end{equation}
This implies that $\con I S q=\con {I\p} S q$. In other words, we can
discard any element $m$ of $\ming I$ where $\proj m$ lies within
$S$. We will apply Equation \eqref{eqn:normalize} whenever it is of
benefit to do so, which it is when $\proj{\ming I}\cap
S\neq\emptyset$. This motivates the following definition.

\begin{definition}[Normal slice]
A slice $(I,S,q)$ is \emph{normal} when $\proj{\ming I}\cap
S=\emptyset$.
\end{definition}

\begin{example}
Let $I$, $p$ and $I_2$ be as in Example \ref{ex:split}. Then
$(I,\ideal p,1)$ is the outer slice after a split on $p$. This slice
is not normal, so we apply Equation \ref{eqn:normalize} to get the
slice $(I_2,\ideal p,1)$, which is the slice $A_2$ from Example
\ref{ex:split}. See Figure $\ref{fig:split}$ for illustrations.
\end{example}

Proposition \ref{prop:split} proves the equations in this section, and
it establishes some results that we will need later.

\begin{proposition}
\label{prop:split}
Let $I$ be a monomial ideal and let $p$ be a monomial. Then
\begin{enumerate}
\item $\gcd(\ming I)$ divides $\gcd(\msm I)$
\label{prop:split:1}

\item $\msm I\cap\ideal p=\msm{I\cap\ideal p}$
\label{prop:split:2}

\item If $p|\gcd(\ming I)$, then $\msm I=\msm{I:p}p$
\label{prop:split:3}

\item $\msm I\cap\ideal p=\msm{I:p}p$
\label{prop:split:4}

\item $\msm I\setminus S=\msm{I\p}\setminus S,\,\,I\p:=\idealBuilder{m\in\ming I}{\proj m\notin S}$
\label{prop:split:5}

\end{enumerate}
\end{proposition}
\begin{proof}
\proofPart{(\ref{prop:split:1})} Let $d\in\msm I$. Let $l_i$ be an $x_i$-label
of $d$ and let $l_j$ be an $x_j$-label of $d$ where $i\neq j$. This is
possible due to the assumption in Section \ref{sec:prelim} that $n\geq 2$. Then
$l_i|dx_i$ and $l_j|dx_j$ so $\gcd(\ming I)|\gcd(l_i,l_j)|d$.

\proofPart{(\ref{prop:split:2})} It follows from Lemma
\ref{lem:splitHelper} below and (\ref{prop:split:1}) that
\[ \msm I\cap\ideal p=\msm{I\cap\ideal p}\cap\ideal p=\msm {I\cap\ideal p} \] 

\proofPart{(\ref{prop:split:3})} If $p|\gcd(\ming I)$ then
$p|\gcd(\msm I)$ by (\ref{prop:split:1}), whereby
\begin{align*}
d\in\msm I
&\biimp
(d/p)p\notin I \text{ and }
(d/p)x_ip\in I\text{ for }i=1,\ldots,n
\\&\biimp
d/p\notin I:p\text{ and }
(d/p)x_i\in I:p\text{ for }i=1,\ldots,n
\\&\biimp
d/p\in\msm{I:p}
\biimp
d\in\msm{I:p}p
\end{align*}

\proofPart{(\ref{prop:split:4})} As $p|\gcd(\ming{I\cap\ideal p})$ and
$(I\cap\ideal p):p=I:p$, we see that
\[
\msm I\cap\ideal p=
\msm{I\cap\ideal p}=
\msm{(I\cap\ideal p):p}p=
\msm{I:p}p \]

\proofPart{(\ref{prop:split:5})} Let $d\in\msm I\setminus S$ and let
$l\in\ming I$ be an $x_i$-label of $d$. Then $l\in\ming{I\p}$ since
$\proj l|d\notin S$. Thus $dx_i\in I\p$ since $l|dx_i$, so
$d\in\msm{I\p}$. Also $d\notin I\supseteq I\p$.

Suppose instead that $d\in\msm{I\p}\setminus S$. Then $dx_i\in
I\p\subseteq I$. If $d\in I$ then there would exist an $m\in\ming
I\setminus\ming {I\p}$ such that $m|d$, which is a contradiction since
then $S\ni\proj m|m|d\notin S$. Thus $d\notin I$ whereby $d\in\msm I$.
\end{proof}

\begin{lemma}
\label{lem:splitHelper}
Let $A$, $B$ and $C$ be monomial ideals. Then $A\cap C=B\cap C$
implies that $\msm A\cap C=\msm B\cap C$.
\end{lemma}
\begin{proof}
Let $d\in\msm A\cap C$. We will prove that $d\in\msm B$.

\proofPart{d\notin B} If $d\in B$ then $d\in B\cap C=A\cap C$ but
$d\notin A$.

\proofPart{d x_i\in B} Follows from $dx_i\in A$ and $d\in C$ since
then $dx_i\in A\cap C=B\cap C$.
\end{proof}

\subsection{The Base Case}
\label{ssec:base}

In this section we present the base case for the Slice Algorithm. A
slice $(I,S,q)$ is a \emph{base case slice} if $I$ is square free or
if $x_1\cdots x_n$ does not divide $\lcm(\ming I)$. Propositions
\ref{prop:baseTrivial} and \ref{prop:base} show why base case slices
are easy to handle.

\begin{proposition}
\label{prop:baseTrivial}
If $x_1\cdots x_n$ does not divide $\lcm(\ming I)$, then
$\msm I=\emptyset$.
\end{proposition}
\begin{proof}
If $\msm I\neq\emptyset$ then there exists some $d\in\msm I$. Let
$m\in\ming I$ be an $x_i$-label of $d$. Then $x_i|m$, so
$x_i|m|\lcm(\ming I)$.
\end{proof}

\begin{proposition}
\label{prop:base}
If $I$ is square free and $I\neq\ideal{x_1,\ldots,x_n}$, then $\msm
I=\emptyset$.
\end{proposition}
\begin{proof}
Let $I$ be square free and let $d\in\msm I$. Let $m_i\in\ming I$ be an
$x_i$-label of $d$ for $i=1,\ldots,n$. Then $d=\proj{\lcm_{i=1}^n
m_i}=1$ so $m_i=x_i$.
\end{proof}

\subsection{Termination And Pivot Selection}
\label{ssec:term}

In this section we show that some quite weak constraints on the choice
of the pivot are sufficient to ensure termination. Thus we leave the
door open for a variety of different pivot selection strategies, which
is something we will have much more to say about in Section
\ref{sec:selection}.

We impose four conditions on the choice of the pivot $p$. These are
presented below, and for each condition we explain why violating that
condition would result in a split that there is no sense in carrying
out. Note that the last two conditions are not necessary at this point
to ensure termination, but they will become so after some of the
improvements in Section \ref{sec:improve} are applied.

\begin{description}
\item[$\boldsymbol{p\notin S}$] If $p\in S$, then the outer slice will
be equal to the current slice.

\item[$\boldsymbol{p\neq 1}$] If $p=1$, then the inner slice will be
equal to the current slice.

\item[$\boldsymbol{p\notin I}$] See Section \ref{ssec:pruneS} and
Equation \eqref{eqn:pruneS} in particular.

\item[$\boldsymbol{p|\proj{\lcm(\ming I)}}$] See Section
\ref{ssec:pruneSMore} and Equation \eqref{eqn:pruneSMore} in
particular.
\end{description}

If a pivot satisfies these four conditions, then we say that it is
\emph{valid}. Proposition \ref{prop:reasonSqFree} shows that non-base
case slices always admit valid pivots, and Proposition \ref{prop:term}
states that selecting valid pivots ensures termination.

\begin{proposition}
\label{prop:reasonSqFree}
Let $(I,S,q)$ be a normal slice for which no valid pivot
exists. Then $I$ is square free.
\end{proposition}
\begin{proof}
Suppose $I$ is not square free. Then there exists an $x_i$ such that
$x_i^2|m$ for some $m\in\ming I$, which implies that $x_i\notin I$.
Also, $x_i\notin S$ since $x_i|\proj m$ and $(I,S,q)$ is normal. We
conclude that $x_i$ is a valid pivot.
\end{proof}

\begin{proposition}
\label{prop:term}
Selecting valid pivots ensures termination.
\end{proposition}
\begin{proof}
Recall that the polynomial ring $R$ is noetherian, so it does not
contain an infinite sequence of strictly increasing ideals. We will
use this to show that the algorithm terminates. Suppose we are
splitting a non-base case slice $A:=(I,S,q)$ on a valid pivot
where $A_1$ is the inner slice and $A_2$ is the outer slice.

Let $f$ and $g$ be functions mapping slices to ideals, and define them
by the expressions $f(I,S,q):=S$ and $g(I,S,q):=\ideal{\lcm(\ming
I)}$. Then the conditions on valid pivots and on non-base case
slices imply that $f(A)\subseteq f(A_1)$, $f(A)\subsetneq f(A_2)$,
$g(A)\subsetneq g(A_1)$ and $g(A)\subseteq g(A_2)$. Also, if we let
$A$ be an arbitrary slice and we let $A\p$ be the corresponding normal
slice, then $f(A)\subseteq f(A\p)$ and $g(A)\subseteq g(A\p)$.

Thus $f$ and $g$ never decrease, and one of them strictly increases on
the outer slice while the other strictly increases on the inner
slice. Thus there does not exist an infinite sequence of splits on
valid pivots.
\end{proof}

\subsection{Pseudo-code}
\label{ssec:pseudo}

This section contains a pseudo-code implementation of the Slice
Algorithm. Note that the improvements in Section \ref{sec:improve} are
necessary to achieve good performance.

The function selectPivot used below returns some valid pivot and can
be implemented according to any of the pivot selection strategies
presented in Section \ref{sec:selection}. A simple idea is to follow
the proof of Proposition \ref{prop:reasonSqFree} and test each
variable $x_1,\ldots,x_n$ for whether it is a valid pivot. If none of
those are valid pivots, then $I\p$ in the pseudo-code below is square
free.

Call the function $\conSym$ below with the parameters $(I, \ideal 0,
1)$ to obtain $\msm I$.

\renewcommand{\labelitemi}{}

\begin{itemize}
\item        function con($I$,\ $S$,\ $q$)
\item \quad let $I\p:=\idealBuilder{m\in\ming I}{\proj m\notin S}$
\item
\item \quad           if $x_1\cdots x_n$ does not divide
                          $\lcm(\ming {I\p})$ then return $\emptyset$
\item \quad           if $I\p$ is square free and 
                          $I\p\neq\ideal{x_1,\ldots,x_n}$ then
                          return $\emptyset$
\item \quad           if $I\p$ is square free and 
                          $I\p=\ideal{x_1,\ldots,x_n}$ then
                          return $\set q$
\item
\item \quad           let $p$ := selectPivot($I\p$,\ $S$)
\item \quad           return $\con {{I\p}:p} {S:p} {qp}\,\,\cup\,\,
                              \con {I\p} {S+\ideal p} q$
\end{itemize}

\label{sec:improve}

\section{Improvements To The Basic Algorithm}

This section contains a number of improvements to the basic version of
the Slice Algorithm presented in Section \ref{sec:basicAlg}.

\subsection{Monomial Lower Bounds On Slice Contents}
\label{ssec:monomialLowerBound}

Let $ql$ be a monomial lower bound on the slice $(I,S,q)$ in the sense
that $ql|d$ for all $d\in\con I S q$. If we then perform a split on
$l$, we can predict that the outer slice will be empty, whereby
Equation \eqref{eqn:split} specializes to Equation
\eqref{eqn:splitLowerBound} below, which shows that we can get the
effect of performing a split while only having to compute a single
slice.
\begin{equation}
\label{eqn:splitLowerBound}
\con I S q = \con {I:l} {S:l} {ql}
\end{equation}

Proposition \ref{prop:split} provides the simple monomial lower bound
$\gcd(\ming I)$, while Proposition \ref{prop:lowerBound} provides a
more sophisticated bound.

\begin{proposition}
\label{prop:lowerBound}
Let $(I,S,q)$ be a slice and let $l(I):=\lcm_{i=1}^n l_i$ where
\[ l_i:=\frac{1}{x_i}\gcd(\ming I\cap\ideal {x_i}) \]
Then $ql(I)$ is a monomial lower bound on $(I,S,q)$.
\end{proposition}
\begin{proof}
Let $d\in\msm I$ and let $m$ be an $x_i$-label of $d$. Then $x_i|m$,
so $l_ix_i|m|dx_i$ whereby $l_i|d$. Thus $l(I)|d$.
\end{proof}

\begin{example}
\label{ex:boundIter}
Let $I:=\ideal{x^2y,xy^2,yz,z^2}$. Then $l(I)=y$ and Equation
\eqref{eqn:splitLowerBound} yields
\[
\con I {\ideal 0} 1=\con {I:y} {\ideal 0} y
\]
where $I:y=\ideal{x^2,xy,z}$. As $l(I:y)=x$ we can apply Equation
\eqref{eqn:splitLowerBound} again to get
\[
\con {I:y} {\ideal 0} y=\con {\ideal{x,y,z}} {\ideal 0}
{xy}=\set{xy}
\]
\end{example}

We can improve on this bound using Lemma \ref{lem:xiMaximal} below.

\begin{definition}[$x_i$-maximal]
A monomial $m\in\ming I$ is \emph{$x_i$-maximal} if
\[ 0<\degx i m=\degx i {\lcm(\ming I)} \]
\end{definition}

\begin{lemma}
\label{lem:xiMaximal}
Let $d\in\msm I$ and let $m$ be an $x_i$-label of $d$. Suppose that
$m$ is $x_j$-maximal for some variable $x_j$. Then $x_i=x_j$.
\end{lemma}
\begin{proof}
Suppose that $x_i\neq x_j$ and let $l$ be an $x_j$-label of $d$. Then
\[ \degx j m\leq \degx j d <
\degx j l\leq\degx j {\lcm(\ming I)}=\degx j m\qedhere\]
\end{proof}

\begin{corollary}
\label{cor:pruneDoubleMaxes}
If $m\in\ming I$ is $x_i$-maximal for two distinct variables, then
$\msm I=\msm{I\p}$ where $I\p:=\ideal{\ming I\setminus\set m}$.
\end{corollary}

\begin{corollary}
\label{cor:improvedLowerBound}
Let $(I,S,q)$ be a slice and let $l_i:= \frac{1}{x_i}\gcd(M_i)$ where
\[ M_i:=\setBuilder{m\in \ming I}
{\text{$x_i$ divides $m$ and $m$ is not $x_j$-maximal for any $x_j\neq
x_i$}} \] Then $q\lcm_{i=1}^n l_i$ is a monomial lower bound on
$(I,S,q)$.
\end{corollary}

It is possible to compute a more exact lower bound by defining
$M_{(i,j)}$ and computing the $\gcd$ of pairs of minimal generators
that could simultaneously be respectively $x_i$ and
$x_j$-labels. However, we expect the added precision to be little and
the computational cost is high. If this is expanded from 2 to $n$
variables, the lower bound is exact, but as costly to compute as the
set $\msm I$ itself.

Corollaries \ref{cor:pruneDoubleMaxes} and
\ref{cor:improvedLowerBound} allow us to make a slice simpler without
changing its content, and they can be iterated until a fixed point is
reached. We call this process \emph{simplification}, and a slice is
\emph{fully simplified} if it is a fixed point of the
process. Proposition \ref{prop:advBaseCase} is an example of how
simplification extends the reach of the base case.

\begin{proposition}
\label{prop:advBaseCase}
Let $A:=(I,S,q)$ be a fully simplified slice. If $\card{\ming I}\leq
n$ then $A$ is a base case slice.
\end{proposition}
\begin{proof}
Assume that $x_1\cdots x_n|\lcm(\ming I)$. Then for each variable
$x_i$, there must be some $m_i\in\ming I$ that is $x_i$-maximal, and
these $m_i$ are all distinct. Since $\card{\ming I}\leq n$ this
implies that $\ming I=\set{m_1,\ldots,m_n}$. Thus $M_i=\set{m_i}$
where $M_i$ is defined in Corollary
\ref{cor:improvedLowerBound}. Furthermore, since $A$ is fully
simplified, $\frac{1}{x_i}\gcd(M_i)=1$, so $m_i=x_i$ and we are done.
\end{proof}

An argument much like that in the proof of Proposition
\ref{prop:advBaseCase} shows that $(I,S,q)$ is a base case if all
elements of $\ming I$ are maximal. If there is exactly one element $m$
of $\ming I$ that is not maximal, then one can construct a new base
case for the algorithm by trying out the possibility of that generator
being an $x_i$-label for each $x_i|m$. One can do the same if there
are $k$ non-maximal elements for any $k\in\N$, but the time complexity
of this is exponential in $k$, so it is slow for large $k$.

Our implementation does this for $k=1,2$, and implementing $k=2$ did
make our program a bit faster. We expect the effect of implementing
$k=3$ would be very small or even negative.

\subsection{Independence Splits}
\label{ssec:independenceSplit}

In this section we define $I$-independence and we show how this
independence allows us to perform a more efficient kind of split. The
content of this section was inspired by a similar technique for
computing Hilbert-Poincar\'e series that was first suggested in
\cite{hseries} and described in more detail in
\cite{bigattiEtAlHilbSerlg}.

\begin{definition}
Let $A,B$ be non-empty disjoint sets such that $A\cup
B=\set{x_1,\ldots,x_n}$.  Then $A$ and $B$ are \emph{$I$-independent}
if $\ming I\cap\ideal A\cap\ideal B=\emptyset$.
\end{definition}

In other words, $A$ and $B$ are $I$-independent if no element of
$\ming I$ is divisible by both a variable in $A$ and a variable in
$B$.

\begin{example}
\label{ex:independenceSplit}
Let $I:=\ideal{x^4,x^2y^2,y^3, z^2,zt,t^2}$. Then $\set{x,y}$ and
$\set{z,t}$ are $I$-independent. It then turns out that we can compute
$\msm I$ independently for $\set{x,y}$ and $\set{z,t}$, which is
reflected in the following equation.
\begin{align*}
\msm I
&=\set{x^3yz,x^3yt,xy^2z,xy^2t}=\set{x^3y,xy^2}\cdot\set{z,t}\\
&=\msm{I\cap\kappa[x,y]}\cdot\msm{I\cap\kappa[z,t]}
\end{align*}
\end{example}

Proposition \ref{prop:independenceSplit} generalizes the observation
in Example \ref{ex:independenceSplit}. The process of applying
Proposition \ref{prop:independenceSplit} is called an
\emph{independence split}.

\begin{proposition}
\label{prop:independenceSplit}
If $A,B$ are $I$-independent, then
\[ \msm I = \msm{I\cap\kappa[A]}\cdot\msm{I\cap\kappa[B]} \]
\end{proposition}
\begin{proof}
Let $A\p:=I\cap\kappa[A]$ and $B\p:=I\cap\kappa[B]$. If $A\p=\ideal 0$
then $\msm I=\emptyset$ by Proposition \ref{prop:baseTrivial}, so we
can assume that $A\p\neq\ideal 0$ and $B\p\neq\ideal 0$. It holds that
\[
\ming I=\ming{A\p}\cup\ming{B\p}
\]
so for monomials $a\in\kappa[A]$ and $b\in\kappa[B]$ we get that
\[
ab\in I\biimp a\in A\p\text{ or }b\in B\p
\]
and thereby
\[
ab\notin I\biimp a\notin A\p\text{ and }b\notin B\p
\]
which implies that
\begin{align*}
ab\in\msm I
\biimp&\,
ab\notin I\text{ and }abx_i\in I\text{ for }x_i\in A\cup B
\\\biimp&\,
a\notin A\p\text{ and }ax_i\in A\p\text{ for }x_i\in A\text{ and}
\\&\,
b\notin B\p\text{ and }bx_i\in B\p\text{ for }x_i\in B
\\\biimp&\,
a\in\msm{A\p}\text{ and }b\in\msm{B\p}\qedhere
\end{align*}
\end{proof}

Given a slice $(I,S,q)$, this brings up the problem of what to do
about $S$ when $A$ and $B$ are $I$-independent but not
$S$-independent. One solution is to remove the elements of $\ming
S\cap\ideal A\cap\ideal B$ from $\ming S$ when doing the independence
split, and then afterwards to remove those computed maximal standard
monomials that lie within $\ideal{\ming S\cap\ideal A\cap\ideal B}$.
Note that this problem does not appear if we use a pivot selection
strategy that only selects pivots of the form $x_i^t$ for $t\in\N$.

\begin{example}
Let $I$ be as in Example \ref{ex:independenceSplit} and consider the
slice $(I,\ideal{x^3y,y^2z},1)$. Then $y^2z$ belongs to neither
$\kappa[x,y]$ nor $\kappa[z,t]$, but we can do the independence split
on the slice $(I,\ideal{x^3y},1)$ which has content
$\set{xy^2}\cdot\set{z,t}=\set{xy^2z,xy^2t}$. We then remove
$\ideal{y^2z}$ from this set, whereby $\con I {\ideal{x^3y,x^2z}}
1=\set{xy^2t}$.
\end{example}

This idea can be improved by observing that when we know $\con {A\p}
{S_A} {q_A}$, we can easily get the monomial lower bound $\gcd(\con
{A\p} {S_A} {q_A})$, and we can exploit this using the technique from
Section \ref{ssec:monomialLowerBound}. This might decrease the size of
$\ming S\cap\ideal A\cap\ideal B$, which can help us compute $\con
{B\p} {S_B} {q_B}$.

This leaves the question of how to detect $I$-independence. This can
be done in space $O(n)$ and nearly in time $O(n\card{\ming I})$ using
the classical union-find algorithm \cite{equivalence,
equivalence2}.\footnote{It can also be done in space $O(n^2)$ and in
time $O(n\card{\ming I}+n^2)$ by constructing a graph in a similar way
and then finding connected components.} See the pseudo-code below,
where $D$ represents a disjoint-set data structure such that
union$(D,x_i,x_j)$ merges the set containing $x_i$ with the set
containing $x_j$. At the end $D$ is the set of independent sets where
$D=\set{\set{x_1,\dots,x_n}}$ implies that there are no independent
sets. The running time claimed above is achieved by using a suitable
data structure for $D$ along with an efficient implementation of
union. See \cite{equivalence,equivalence2} for details.

\begin{itemize}
\item \quad           let $D:=\set{\set{x_1},\ldots,\set{x_n}}$.
\item \quad           for each $m\in\ming I$ do
\item \quad\quad        pick an arbitrary $x_i$ that divides $m$
\item \quad\quad        for each $x_j$ that divides $m$ do
\item \quad\quad\quad     union$(D,x_i,x_j)$
\end{itemize}

This is an improvement on the $O(n^2\card{\ming I})$ algorithm for
detecting independence suggested in \cite{bigattiEtAlHilbSerlg}. That
algorithm is similar to the one described here, the main difference
being the choice of data structure.

\subsection{A Base Case Of Two Variables}

When $n=2$ there is a well known and more efficient way to compute
$\msm I$. This is also useful when an independence split has reduced
$n$ down to two.

Let $\set{m_1,\ldots,m_k}:=\ming I$ where $m_1$,\ldots,$m_k$ are
sorted in ascending lexicographic order where $x_1>x_2$. Let
$\tau(x^u,x^v):=x^{(v_1,u_2)}$. Then
\[ \msm I=\set{\tau(m_1,m_2),\tau(m_2,m_3),\ldots,\tau(m_{k-1},m_k)}
. \]

\subsection{Prune $S$}
\label{ssec:pruneS}

Depending on the selection strategy used, it is possible for the $S$
in $(I,S,q)$ to pick up a large number of minimal generators, which
can slow things down. Thus there is a point to removing elements of
$\ming S$ when that is possible without changing the content of the
slice. Equation \eqref{eqn:pruneS} allows us to do this.
\begin{equation}
\label{eqn:pruneS}
\con I S q=\con I {S\p} q,\quad S\p:=\idealBuilder{m\in\ming S}{m\notin I}
\end{equation}

\begin{example}
Consider the slice $(\ideal{x^2,y^2,z^2,yz},\ideal{xyz},1)$. Then
$p:=x$ is a valid pivot, yielding the inner slice
$(\ideal{x,y^2,z^2,yz},\ideal{yz},x)$. We can now apply Equation
\eqref{eqn:pruneS} to turn this into $(\ideal{x,y^2,z^2,yz},\ideal{0},x)$.
\end{example}

Proposition \ref{prop:term} states that the Slice Algorithm
terminates, and we need to prove that this is still true when we use
Equation \eqref{eqn:pruneS}. Fortunately, the same proof can be used,
except that the definition of the function $f$ needs to be changed
from $f(I,S,q)=S$ to $f(I,S,q):=I+S$. Note that the condition on a
valid pivot $p$ that $p\notin I$ is there to make this work.

\subsection{More Pruning of $S$} 
\label{ssec:pruneSMore}

We can prune $S$ using Equation \eqref{eqn:pruneSMore} below, and for
certain splitting strategies this will even allow us to never add
anything to $S$.
\begin{equation}
\label{eqn:pruneSMore}
\con I S q =\con I {S\p} q, S\p := \idealBuilder{m\in\ming S}{m\text{ divides }\proj{\lcm(\ming I)}}
\end{equation}
To prove this, observe that any $d\in\con I S q$ divides
$\proj{\lcm(\ming I)}$.

\begin{example}
\label{ex:pruneSMore}
Consider the slice $(\ideal{x^2,xy,y^2},\ideal 0,1)$. Then $p:=x$ is a
valid pivot, yielding the normalized outer slice
$(\ideal{xy,y^2},\ideal x,1)$. We can now apply Equation
\eqref{eqn:pruneSMore} to turn this into $(\ideal{xy,y^2},\ideal
0,1)$.
\end{example}

Similarly, Equation \eqref{eqn:pruneSMore} will remove any generator
of the form $x_i^t$ from $S$. So if we use a pivot of the special form
$p=x_i^t$, and we apply a normalization and Equation
\eqref{eqn:pruneSMore} to the outer slice, we can turn Equation
\eqref{eqn:split} into
\[
\con I S q=\con {I:p} {S:p} {qp}\,\,\cup\,\,
\con {\ideal{\ming I\setminus\ideal{px_i}}} S q
\]
which for $S=\ideal 0$ and $q=1$ specializes to
\[
\msm I=\msm{I:x_i^t}x_i^t\,\,\cup\,\,
\msm{\ideal{\ming I\setminus\ideal{x_i^{t+1}}}}
\]
An implementer who does not want to deal with $S$ might prefer this
equation to the more general Equation \eqref{eqn:split}.

We need to prove that the algorithm still terminates when using
equations \eqref{eqn:pruneS} and \eqref{eqn:pruneSMore}. We can use
the same proof as in Proposition \ref{prop:term}, except that we need
to replace the definition of $f$ from that proof with
$f(I,S,q):=I+S+\ideal{x_1^{u_1},\ldots,x_n^{u_n}}$ where
$x^u:=\lcm(\ming I)$. Note that the condition on a valid pivot $p$
that $p|\proj{\lcm(\ming I)}$ is there to make this work.

\subsection{Minimizing The Inner Slice}
\label{ssec:minimize}

A time-consuming step in the Slice Algorithm is to compute $I:p$ for
each inner slice $(I:p,S:p,qp)$. By \emph{minimizing}, we mean the
process of computing $\ming{I:p}$ from $\ming I$, which is done by
removing the non-minimal elements of $\ming I : p :=
\setBuilder{m:p}{m\in\ming I}$ where $m:p:=\frac{m}{\gcd(m,p)}$.

Proposition \ref{prop:minimize} below makes it possible to do this
using fewer divisibility tests than would otherwise be required. As
seen by Corollary \ref{cor:minimize} below, this generalizes both
statements of \cite[Proposition 1]{bigattiHilbSerComp} from $p$ of the
form $x_i^t$ to general $p$.\footnote{This provides an answer to the
statement from \cite[p. 11]{bigattiHilbSerComp} that ``These remarks
drastically reduce the number of divisibility tests, but they do not
easily generalize for non-simple-power pivots, not even for
power-products with only two indeterminates.''} See \cite[Section
6]{bigattiEtAlHilbSerlg} for an even earlier form of these ideas.

Note that the techniques in this section also apply to computing
intersections $I\cap\ideal p$ of a momomial ideal with a principal
ideal generated by a monomial, since $\ming{I\cap\ideal
p}=\setBuilder{\lcm(m,p)}{m\in\ming{I:p}}$.

The most straightforward way to minimize $\ming I:p$ is to consider
all pairs of distinct $a,b\in\ming I:p$ and then to remove $b$ if
$a|b$. It is well known that this can be improved by sorting $\ming
I:p$ according to some term order, in which case a pair only needs to
be considered if the first term comes before the last. This halves the
number of divisibility tests that need to be carried out.

We can go further than this, however, because we know that $\ming I$
is already minimized. Proposition \ref{prop:minimize} shows how we can
make use of this information.

\begin{proposition}
\label{prop:minimize}
Let $x^a$, $x^b$ and $x^p$ be monomials such that $x^a$ does not
divide $x^b$. Then $x^a:x^p$ does not divide $x^b:x^p$ if it holds for
$i=1,\ldots,n$ that $p_i<a_i\lor a_i\leq b_i$.
\end{proposition}
\begin{proof}
We prove the contrapositive statement, so suppose that $x^a$ does not
divide $x^b$ and that $x^u:=x^a:x^p$ divides $x^v:=x^b:x^p$. Then
there is an $i$ such that $a_i>b_i$. As $\max(p_i,a_i)=u_i+p_i\leq
v_i+p_i=\max(p_i,b_i)$ we conclude that $p_i\geq a_i$.
\end{proof}

This allows us to draw some simple and useful conclusions.

\begin{corollary}
\label{cor:minimize}
Let $a,b\in\ming I$. Then $a:p$ does not divide $b:p$ if any one of
the following two conditions is satisfied.
\begin{enumerate}
\item $\sqrt a = \sqrt{a:p}$
\label{xprop:minimize:1}

\item $\gcd(a,p)|\gcd(b,p)$
\label{xprop:minimize:2}

\end{enumerate}
\end{corollary}

\begin{corollary}
\label{cor:minimizeSupp}
If $a\in\ming I$ and $p|\proj a$, then $a:p\in\ming{I:p}$ and $a:p$
does not divide any other element of $\ming I:p$.
\end{corollary}

We can push Proposition \ref{prop:minimize} further than this. Fix
some monomial $p$, let $x^t:=p$, $t\in\N^n$, and define the function
$f$, which maps monomials to vectors in $\N^n$, by
\[
(f(x^u))_i=
\begin{cases}
0,&\textrm{for }t_i=0, \\
\min(u_i,t_i+1),&\textrm{for }t_i\neq 0.
\end{cases}
\]
Also, define the relation $u\prec v$ for vectors $u,v\in\N^n$ by
\[
u\prec v\text{ if there exists an $i\in\set{1,\ldots,n}$ such that
$t_i\geq u_i>v_i$.}
\]
By Proposition \ref{prop:minimize}, this implies that if $a,b\in\ming
I$ and $f(a)\nprec f(b)$, then $a:p$ does not divide $b:p$. This
motivates us to define the function $L(a)$ by
\[
L(a):=\setBuilder{m\in\ming I}{f(a)=f(m)}.
\]
If $a,b\in\ming I$ and $f(a)\nprec f(b)$, then no element of $L(a):p$
divides any element of $L(b):p$. In particular, no element of $L(a):p$
divides any other. Note that the domain of $L$ is a partition of
$\ming I$.

This technique works best when most of the non-empty sets $L(a)$
contain considerably more than a single element, which is likely to be
true e.g. if $p$ is a small power of a single variable. Even in cases
where most of the non-empty sets $L(a)$ consist of only a few
elements, it will likely still pay off to consider $L(1)$ and to make
use of Corollary \ref{cor:minimizeSupp}.

\begin{example}
Let $I:=\ideal{x^5y,x^2y^2,x^2z^3,xy^3,xyz^3,yz^2}$ and $p:=x^3$.
Then
\begin{align*}
L(x^4)&=\set{x^5y}          &L(x^4):p&=\set{x^2y}\\
L(x^2)&=\set{x^2y^2,x^2z^3} &L(x^2):p&=\set{y^2,z^3}\\
L(x)  &=\set{xy^3,xyz^3}    &L(x):p  &=\set{y^3,yz^3}\\
L(1)  &=\set{yz^2}          &L(1):p  &=\set{yz^2}
\end{align*}
We will process these sets from the top down. The set $L(x^4)$ is
easy, since $p|\proj{x^5y}$, so we do not have to do any divisibility
tests for $x^5y$.

Then comes $L(x^2)$. We have to test if any elements of $L(x^2):p$
divide any elements of $L(x):p$ or $L(1):p$. It turns out that
$x^2y^2:p|xy^3:p$ and $x^2z^3:p|xyz^3:p$, so we can remove all of
$L(x)$ from consideration.  We do not need to do anything more for
$L(1)$, so we conclude that $\ming{I:p}=\ideal{x^2y,y^2,z^3,yz^2}$.
\end{example}

\subsection{Reduce The Size Of Exponents}
\label{ssec:remap}

Some applications require the irreducible decomposition of monomial
ideals $I$ where the exponents that appear in $\ming I$ are very
large. One example of this is the computation of Frobenius numbers
\cite{bfrob, frobPoint}.

This presents the practical problem that these numbers are larger than
can be natively represented on a modern computer. This necessitates
the use of an arbitrary precision integer library, which imposes a
hefty overhead in terms of time and space. One solution to this
problem is to report an error if the exponents are too large, as
indeed the programs Monos \cite{monos} and Macaulay 2 \cite{M2} do for
exponents larger than $2^{15}-1$ and $2^{31}-1$ respectively.

In this section, we will briefly describe how to support arbitrarily
large exponents without imposing any overhead except for a quick
preprocessing step. The most time-consuming part of this preprocessing
step is to sort the exponents.

Let $f$ be a function mapping monomials to monomials such that
$f(ab)=f(a)f(b)$ when $\gcd(a,b)=1$. Suppose that $a|b\Rightarrow
f(a)|f(b)$ and that $f$ is injective for each $i$ when restricted to
the set $\setBuilder{x_i^{v_i}}{x^v\in\ming I}$. The reader may verify
that then
\[
x_1\cdots x_n\msm I=f^{-1}(x_1\cdots x_n\msm{\ideal{f(\ming I)}})
.\]
The idea is to choose $f$ such that the exponents in $f(\ming I)$ are
as small as possible, which can be done by sorting the exponents that
appear in $\ming I$. If this is done individually for each variable,
then $\card{\ming I}$ is the largest integer that can appear as an
exponent in $f(\ming I)$. Thus we can compute $\msm I$ in terms of
$\msm{\ideal{f(\ming I)}}$, which does not require large integer
computations.

\begin{example}
If $I:=\ideal{x^{100},x^{40}y^{20},y^{90}}$ then we can choose the
function $f$ such that $\ideal{f(\ming I)}=\ideal{x^2,xy,y^2}$.
\end{example}

The underlying mathematical idea used here is that it is the order
rather than the value of the exponents that matters. This idea can
also be found e.g. in \cite[Remark 4.6]{scarfComplex}, though in the
present paper it is used to a different purpose.

\subsection{Label Splits}

In this section we introduce \emph{label splits}. These are based on
some properties of labels which pivot splits do not make use of.

Let $(I,S,q)$ be the current slice, and assume that it is fully
simplified and not a base case slice. The first step of a label split
is then to choose some variable $x_i$ such that $\ming I\cap\ideal{x_i}\neq\set{x_i}$. Let $L:=\setBuilder{x^u\in\ming I}{u_i=1}$. Then $L$ is non-empty
since the current slice is fully simplified. Assume for now that
$\card L=1$ and let $l\in L$.

Observe that if $d\in\msm I$, then $\frac{l}{x_i}|d$ if and only if
$l$ is an $x_i$-label of $d$, which is true if and only if $x_i$ does
not divide $d$. This and Equation \eqref{eqn:sliceIbyM} implies that
\[
\con I S 1\setminus\ideal{x_i}=
\con I S 1\cap\ideal{\frac{l}{x_i}}=
\con{I:\frac{l}{x_i}} {S:\frac{l}{x_i}} {\frac{l}{x_i}}
\]
\[
\con I S 1\cap\ideal{x_i}=\con {I:x_i} {S:x_i} {x_i}
\]
whereby
\[
\con I S q=\con {I:x_i} {S:x_i} {qx_i}\,\,\cup\,\,
\con{I:\frac{l}{x_i}} {S:\frac{l}{x_i}} {q\frac{l}{x_i}}
\]
This equation describes a label split on $x_i$ in the case where
$\card L=1$. In general $\card L$ can be larger than one, so let
$L=\set{l_1,\ldots,l_k}$ and define
\[ I_j:=I\colon{\frac{l_j}{x_i}},\quad
S_j:=\left(S+\ideal{\frac{l_1}{x_i},\ldots,\frac{l_{j-1}}{x_i}}\right):{\frac{l_j}{x_i}},\quad
q_j:=q{\frac{l_j}{x_i}} \] for $j=1,\ldots,k$. Then $\con {I_j} {S_j}
{q_j}$ is the set of those $d\in\con I S q$ such that $l_j$ is an
$x_i$-label of $d$, and such that none of the monomials
$l_1,\ldots,l_{j-1}$ are $x_i$-labels of $d$. This implies that
\begin{equation}
\label{eqn:labelSplit}
\con I S q=\con {I:x_i} {S:x_i} {qx_i}\,\,\bigcup_{j=1}^k
\con {I_j} {S_j} {q_j}
\end{equation}
where the union is disjoint. This equation defines a label split on
$x_i$.

An advantage of label splits is that if $I$ is artinian, $S=\ideal 0$
and $\card L=1$, then none of the slices on the right hand side of
Equation \eqref{eqn:labelSplit} are empty. These conditions will
remain true throughout the computation if the ideal is artinian and
generic and we perform only label and independence splits. Example
\ref{ex:nongenericLabelSplit} shows that a label split can produce
empty slices when $\card L>1$.

\begin{example}
\label{ex:nongenericLabelSplit}
Let $I:=\ideal{x^4,y^4,z^4,xy,xz}$. We perform a label split on $x$
where $l_1:=xy$ and $l_2:=xz$, which yields the following equation.
\begin{align*}
\con I {\ideal 0} 1
&=
\con {\ideal{x^3,y,z}} {\ideal 0} x   & \text{(this is $(I:x_i, S:x_i, qx_i)$)}
\\&\cup
\con {\ideal{x,y^3,z^4}} {\ideal 0} y & \text{(this is $(I_1,S_1,q_1)$)}
\\&\cup
\con {\ideal{x,y^4,z^3}} {\ideal y} z & \text{(this is $(I_2,S_2,q_2)$)}
\\&=
\set{x^3}\cup\set{y^3z^3}\cup\emptyset=\set{x^3,y^3z^3}
\end{align*}
The reason that $(I_2,S_2,q_2)$ is empty is that both $l_1$ and $l_2$
are $x$-labels of $y^3z^3$.
\end{example}

Using only label splits according to the VarLabel strategy discussed
in Section \ref{sssec:labelSelect} makes the Slice Algorithm behave as
a version of the Label Algorithm \cite{rouneLabelAlg}. See the
External Corner Algorithm \cite{frobPoint} for an earlier form of some
of the ideas behind the Label Algorithm.

\section{Split Selection Strategies}
\label{sec:selection}

We have not specified how to select the pivot monomial when doing a
pivot split, or when to use a label split and on what variable. The
reason for this is that there are many possible ways to do it, and it
is not clear which one is best. Indeed, it may be that one split
selection strategy is far superior to everything else in one
situation, while being far inferior in another. Thus we examine
several different selection strategies in this section.

We are in the fortunate situation that an algorithm for computing
Hilbert-Poincar\'e series has an analogous issue of choosing a pivot
\cite{bigattiEtAlHilbSerlg}. Thus we draw on the literature on that
algorithm to get interesting pivot selection strategies
\cite{bigattiEtAlHilbSerlg, bigattiHilbSerComp}, even though these
strategies do have to be adapted to work with the Slice Algorithm. The
independence and label strategies are the only ones among the
strategies below that is not similar to a known strategy for the
Hilbert-Poincar\'e series algorithm.

It is assumed in the discussion below that the current slice is fully
simplified and not a base case slice. Note that all the strategies
select valid pivots only. We examine the practical merit of these
strategies in Section \ref{ssec:strategyBench}.

\subsection{The Minimal Generator Strategy}

We abbreviate this as \emph{MinGen}.

\subsubsection*{Selection}

This strategy picks some element $m\in\ming I$ that is not square free
and then selects the pivot $\proj m$.

\subsubsection*{Analysis}

This strategy chooses a pivot that is maximal with respect to the
property that it removes at least one minimal generator from the outer
slice. This means that the inner slice is easy, while the outer slice
is comparatively hard since we can be removing as little as a single
minimal generator.

\subsection{The Pure Power Strategies} There are three pure power strategies.
\label{sssec:purePowerSelect}

\subsubsection*{Selection}
These strategies choose a variable $x_i$ that maximizes $\card{\ming
I\cap\ideal {x_i}}$ provided that $x_i^2|\lcm(\ming I)$. Then they
choose some positive integer $e$ such that $x_i^{e+1}|\lcm(\ming I)$
and select the pivot $x_i^e$.

The strategy \emph{Minimum} selects $e:=1$ and the strategy
\emph{Maximum} selects $e:=\degx i {\lcm(\ming I)}-1$. The strategy
$\emph{Median}$ selects $e$ as the median exponent of $x_i$ from the
set $\ming I\cap\ideal{x_i}$.

Note that the Minimum strategy makes the Slice Algorithm behave as a
version of the staircase-based algorithm due to Gao and Zhu
\cite{artinianStairIrr}.

\subsubsection*{Analysis}

The pure power strategies have the advantage that the minimization
techniques described in Section \ref{ssec:minimize} work especially
well for pure power pivots. Maximum yields an easy inner slice and a
hard outer slice, while Minimum does the opposite. Median achieves a
balance between the two.

\subsection{The Random GCD Strategy}
We abbreviate this as \emph{GCD}.

\subsubsection*{Selection}

Let $x_i$ be a variable that maximizes $\card{\ming I\cap\ideal{x_i}}$
and pick three random monomials $m_1,m_2,m_3\in\ming
I\cap\ideal{x_i}$. Then the pivot is chosen to be
$p:=\proj{\gcd(m_1,m_2,m_3)}$. If $p=1$, then the GCD strategy fails,
and we might try again or use a different selection strategy.

\subsubsection*{Analysis}

We consider this strategy because a similar strategy has been found to
work well for the Hilbert-Poincar\'e series algorithm mentioned above.

\subsection{The Independence Strategy}
We abbreviate this as \emph{Indep}.

\subsubsection*{Selection}

The independence strategy picks two distinct variables $x_i$ and
$x_j$, and then selects the pivot
$p:=\proj{\gcd(\ming{I}\cap\ideal{x_ix_j})}$. If $p=1$, then the
independence strategy fails, and we might try again or use a different
selection strategy.

\subsubsection*{Analysis}

The pivot $p$ is the maximal monomial that will make every minimal
generator that is divisible by both $x_i$ and $x_j$ disappear from the
outer slice. The idea behind this is to increase the chance that we
can perform an independence split on the outer slice while having a
significant impact on the inner slice as well.

\subsection{The Label Strategies} There are several label strategies.
\label{sssec:labelSelect}

\subsubsection*{Selection}
These strategies choose a variable $x_i$ such that $\ming I\cap
\ideal{x_i}\neq\set{x_i}$ and then perform a label split on $x_i$. The
strategy MaxLabel maximizes $\card{\ming I\cap\ideal{x_i}}$,
\mbox{VarLabel} minimizes $i$ and MinLabel minimizes
$\card{\setBuilder{x^u\in\ming I}{v_i=1}}$ while breaking ties
according to MaxLabel.

Note that the VarLabel strategy makes the Slice Algorithm behave as a
version of the Label Algorithm \cite{rouneLabelAlg}.

\subsubsection*{Analysis}
MaxLabel chooses the variable that will have the biggest impact, while
MinLabel avoids considering as many empty slices by keeping
$\card{\ming S}$ small. MinLabel is being considered due to its
relation to the Label Algorithm.

\section{Applications To Optimization}
\label{sec:apps}

Sometimes we compute a socle or an irreducible decomposition because
we want to know some property of it rather than because we are
interested in knowing the socle or decomposition itself. This kind of
situation often has the form
\[
\text{maximize }v(J)\text{ subject to }J\in\decom I
\]
where $v$ is some function mapping $\decom I$ to $\R$. We call such a
problem an Irreducible Decomposition Program (IDP). As described in
sections \ref{ssec:linearIDP} and \ref{ssec:igap}, applications of IDP
include computing the integer programming gap, Frobenius numbers and
the codimension of a monomial ideal.

The Slice Algorithm can solve some IDPs in much less time than it
would need to compute all of $\decom I$, and that is the subject of
this section. Section \ref{ssec:branch} explains the general principle
of how to do this, while Section \ref{ssec:monomialBounds} provides
some useful techniques for making use of the principle. Sections
\ref{ssec:linearIDP} and \ref{ssec:igap} present examples of how to
apply these techniques.

\subsection{Branch And Bound Using The Slice Algorithm}
\label{ssec:branch}

In this section we explain the general principle of solving IDPs using
the Slice Algorithm.

The first issue is that the Slice Algorithm is concerned with
computing maximal standard monomials while IDPs are about irreducible
decomposition. We deal with this by using the function $\phi$ from
Section \ref{ssec:msm} to reformulate an IDP of the form
\[
\text{maximize }v\p(J)\text{ subject to }J\in\decom {I\p}
\]
into the form
\[
\text{maximize }v(d)\text{ subject to }d\in\msm {I}
\]
where $v(d):=v\p(\phi(d))$ and $I:=I\p+\ideal{x_1^t,\ldots,x_n^t}$ for
some $t>>0$.

It is a simple observation that there is no reason to compute all of
$\msm I$ before beginning to pick out the element that yields the
greatest value of $v$. We might as well not store $\msm I$, and only
keep track of the greatest value of $v$ found so far.

We define a function $b(I,S,q)$ that maps slices $(I,S,q)$ to real
numbers to be an \emph{upper bound} if $d\in\con I S q$ implies that
$v(d)\leq b(I,S,q)$. We will now show how to use such an upper bound
$b$ to turn the Slice Algorithm into a \emph{branch and bound}
algorithm.

Suppose that the Slice Algorithm is computing the content of a slice
$(I,S,q)$, and that $b(I,S,q)$ is less than or equal to the greatest
value of $v$ found so far. Then we can skip the computation of $\con I
S q$, since no element of $\con I S q$ improves upon the greatest
value of $v$ found so far.

We can take this a step further by extending the idea of monomial
lower bounds from Section \ref{ssec:monomialLowerBound}. The point
there was that if we can predict that the outer slice of some pivot
split will be empty, then we should perform that split and ignore the
outer slice. That way we get the benefit of a split while only having
to examine a single slice. In the same way, if we can predict that one
slice of some pivot split will not be able to improve upon the best
value found so far, we should perform the split and ignore the
non-improving slice. The hard part is to come up with a way to find
pivots where such a prediction can be made. Sections
\ref{ssec:linearIDP} and \ref{ssec:igap} provide examples of how this
can be done.

A prerequisite for applying the ideas in this section is to construct
a bound $b$. It is not possible to say how to do this in general,
since it depends on the particulars of the problem at hand, but
Section \ref{ssec:monomialBounds} presents some ideas that can be
helpful.

\subsection{Monomial Bounds}
\label{ssec:monomialBounds}

In this section we present some ideas that can be useful when
constructing upper bounds for IDPs of the form
\[
\text{maximize }v(d)\text{ subject to }d\in\msm I
.\]

Suppose that $v$ is decreasing in the sense that if $a|b$ then
$v(a)\geq v(b)$. Then $b(I,S,q):=v(q)$ is an upper bound, since if
$d\in\con I S q$ then $q|d$ so $v(d)\leq v(q)$.

Suppose instead that $v$ is increasing in the sense that if $a|b$ then
$v(a)\leq v(b)$. Then $b(I,S,q):=v(q\proj{\lcm(\ming I)})$ is an upper
bound, since if $d\in\con I S q$ then $d|q\proj{\lcm(\ming I)}$ by
Proposition \ref{prop:monomialUpperBound} below, so $v(d)\leq
v(q\proj{\lcm(\ming I)})$. Any monomial upper bound on $\con I S q$
yields an upper bound in the same way.

\begin{proposition}
\label{prop:monomialUpperBound}
If $d\in\msm I$ then $d|\proj{\lcm(\ming I)}$.
\end{proposition}
\begin{proof}
Let $d\in\msm I$ and let $m_i\in\ming I$ be an $x_i$-label of $d$ for
$i=1,\ldots,n$. Then $d=\proj{\lcm_{i=1}^n m_i}$ divides
$\proj{\lcm(\ming I)}$.
\end{proof}

Sections \ref{ssec:linearIDP} and \ref{ssec:igap} provide examples of
how these ideas can be applied.

\subsection{Linear IDPs, Codimension And Frobenius Numbers}
\label{ssec:linearIDP}

Let $r\in\R^n$ and define the function $v_r(x^u):=u\cdot r$. Then we
refer to IDPs of the form \eqref{eqn:linearIDP} as \emph{linear}.
\begin{equation}
\label{eqn:linearIDP}
\text{maximize }v_r(d)\text{ subject to }d\in\msm I
\end{equation}

It is well known that the \emph{codimension} of a monomial ideal $I\p$
equals the minimal number of generators of the ideals in $\decom
{I\p}$. The reader may verify that this is exactly the optimal value
of the IDP \eqref{eqn:linearIDP} if we let
$I:=\sqrt{I\p}+\ideal{x_1^2,\ldots,x_n^2}$ and $r=(1,\ldots,1)$,
noting the well known fact that the codimension of an ideal does not
change by taking the radical. This implies that solving IDPs is
NP-hard since computing codimensions of monomial ideals is NP-hard
\cite[Proposition 2.9]{hseries}. Linear IDPs are also involved in the
computation of Frobenius numbers \cite{bfrob, frobPoint}.

Let us return to the general situation of $r$ and $I$ being
arbitrary. Our goal in this section is to solve IDPs of the form
\eqref{eqn:linearIDP} efficiently by constructing a bound. The
techniques from Section \ref{ssec:monomialBounds} do not immediately
seem to apply, since $v_r$ need neither be increasing nor
decreasing. To deal with this problem, we will momentarily restrict
our attention to some special cases.

Let $a\in\R^n_{\geq 0}$ be a vector of $n$ non-negative real numbers,
and define $v_a(x^u):=u\cdot a$. We will construct a bound for the IDP
\[
\text{maximize }v_a(d)\text{ subject to }d\in\msm I
.\]

This is now easy to do, since $v_a$ is increasing so that we can use
the techniques from Section \ref{ssec:monomialBounds}. Specifically,
$v_a(d)\leq v_a(q\proj{\lcm(\ming I)})$ for all $d\in\con I S q$.

Similarly, let $b\in\R^n_{\leq0}$ be a vector of $n$ non-positive real
numbers, and define $v_b(x^u):=u\cdot b$. We will construct a bound for the IDP
\[
\text{maximize }v_b(d)\text{ subject to }d\in\msm I
.\]
This is also easy, since $v_b$ is decreasing so that we can use the
techniques from Section \ref{ssec:monomialBounds}. Specifically,
$v_b(d)\leq v_b(q)$ for all $d\in\con I S q$.

We now return to the issue of constructing a bound for the IDP
\eqref{eqn:linearIDP}. Choose $a\in\R^n_{\geq0}$ and
$b\in\R^n_{\leq0}$ such that $r=a+b$. Then we can combine the bounds
for $v_a$ and $v_b$ above to get a bound for $v$. So if $d\in\con I S
q$, then
\[
v(d)=v_a(d)+v_b(d)\leq v_a(q\proj{\lcm(\ming I)})+v_b(q) =: b(I,S,q)
\]

Now that we have a bound $b$, we follow the suggestion from Section
\ref{ssec:branch} that we should devise a way to find pivots where we
can predict that one of the slices will be non-improving. Let
$(I,S,q)$ be the current slice and let $x^u:=\lcm(\ming I)$.

Suppose that $r_i$ is positive and consider the outer slice
$(I\p,S\p,q\p)$ from a pivot split on $x_i$. We can predict that the
exponent of $x_i$ in our monomial upper bound will decrease from
$\degx i q+u_i-1$ down to $\degx i q$. Thus we get that
\[
r_i(u_i - 1)\leq b(I,S,q)-b(I\p,S\p,q\p)
,\]
whereby
\[
b(I\p,S\p,q\p)\leq b(I,S,q)-r_i(u_i - 1)
,\]
which implies that the outer slice is non-improving if
\begin{equation}
\label{eq:outerBound}
b(I,S,q)-r_i(u_i - 1) \leq \tau
,\end{equation}
where $\tau$ is the best value found so far. We can do a similar thing
if $r_i$ is negative by considering the value of $\degx i {q\p}$ on
the inner slice of a pivot split on $x_i^{u_i-1}$.

As we will see in Section \ref{ssec:boundBench}, this turns out to
make things considerably faster. One reason is that checking Equation
\eqref{eq:outerBound} for each variable $x_i$ is very fast, because it
only involves computations on the single monomial $\lcm(\ming
I)$. Another reason is that we can iterate this idea, as moving to the
inner or outer slice can reduce the bound, opening up the possibility
for doing the same thing again. We can also apply the simplification
techniques from Section \ref{ssec:monomialLowerBound} after each
successful application of Equation \eqref{eq:outerBound}.

\subsection{The Integer Programming Gap}
\label{ssec:igap}

Let $c\in\Q^n$ and $d\in\Z^k$, and let $A$ be a $k\times n$ integer
matrix. The \emph{integer programming gap} of a bounded and feasible
integer program of the form
\[
\text{minimize }c\cdot x\text{ subject to }Ax=d,\ x\in\N^n
\]
is the difference between its optimal value and the optimal value of
its \emph{linear programming relaxation}, which is defined as the
linear program
\[
\text{minimize }c\cdot x\text{ subject to }Ax=d,\ x\in\R^n_{\geq 0}
.\]

The paper \cite{igap} describes a way to compute the integer
programming gap that involves the sub-step of computing an irreducible
decomposition $\decom{I\p}$ of a monomial ideal $I\p$. Our goal in
this section is to show that this sub-step can be reformulated as an
IDP whose objective function $v$ satisfies the property that $a|b\imp
v(a)\leq v(b)$ whereby we can construct a bound using the technique
from Section \ref{ssec:monomialBounds}.

First choose $t>>0$ and let $I:=I\p+\ideal{x_1^{t+1},\ldots,x_n^{t+1}}$ so
that we can consider $\msm I$ in place of $\decom{I\p}$. Define
$\psi\colon \N^n\mapsto \N^n$ by the expression
\[
(\psi(u))_i :=
\begin{cases}
u_i,&\textrm{for }u_i < t, \\
0,&\textrm{for }u_i\geq t.
\end{cases}
\]
So if $t=4$ then $\psi(3,4,5)=(3,0,0)$. Define $v(u)$ for $u\in\N^n$
as the optimal value of the following linear program. We say that this
linear program is \emph{associated to} $u$.
\[
\begin{array}{rl}
\text{maximize }&c\cdot(\psi(u) - w)\\
\text{subject to }&A(\psi(u)-w)=0,\ w\in\R^n
\\
\text{and }&w_i\geq 0\text{ for those }i\text{ where }u_i<t 
\end{array}
\]

The IDP that the algorithm from \cite{igap} needs to solve is then
\[ \text{maximize }v(u)\text{ subject to }x^u\in\msm I. \]
By Proposition \ref{prop:niceGap} below, we can construct a bound for
this IDP using the technique from Section
\ref{ssec:monomialBounds}. Note that we can use this bound to search
for non-improving outer slices for pivots of the form $x_i$ in the
exact same way as described for linear IDPs in Section
\ref{ssec:linearIDP}.

\begin{proposition}
\label{prop:niceGap}
The function $v$ satisfies the condition that $x^a|x^b\imp v(a)\leq
v(b)$.
\end{proposition}
\begin{proof}
Let $e_i\in\N^n$ be a vector of zeroes except that the $i$'th entry is
1. It suffices to prove that $v(u)\leq v(u+e_i)$ for $u\in\N^n$.  Let
$w\in\R^n$ be some optimal solution to the linear program associated
to $u$. We will construct a feasible solution $w\p$ to the linear
program associated to $u+e_i$ that has the same value. We will ensure
this by making $w\p$ satisfy the equation $\psi(u)-w=\psi(u+e_i)-w\p$.

\proofPart{\text{The case }u_i+1< t} Let $w\p:=w+e_i$.

\proofPart{\text{The case }u_i+1=t} Let $w\p:=w-u_ie_i$. Note that
the non-negativity constraint on the $i$'th entry of $w\p$ is lifted
due to $u_i+1=t$.

\proofPart{\text{The case }u_i+1> t} Let $w\p:=w$. Note that
this case is not relevant to the computation since no upper bound will
be divisible by $x_i^{t+1}$.
\end{proof}

\section{Benchmarks}
\label{sec:bench}

We have implemented the Slice Algorithm in the software system Frobby
\cite{frobby}, and in this section we use Frobby to explore the Slice
Algorithm's practical performance. Section \ref{ssec:dataBench}
describes the test data we use, Section \ref{ssec:strategyBench}
compares a number of split selection strategies, Section
\ref{ssec:otherBench} compares Frobby to other programs and finally
Section \ref{ssec:boundBench} evaluates the impact of the bound
optimization from Section \ref{sec:apps}.

\subsection{The Test Data}
\label{ssec:dataBench}

In this section we briefly describe the test data that we use for the
benchmarks. Table \ref{fig:inputData} displays some standard
information about each input. The data used is publicly available at
\url{http://www.broune.com/}.

\begin{table}
\centering
\begin{tabular}{|l||r|r|r|r|}
\hline
name & n & $|\min(I)|$ & $|\textrm{irr}(I)|$ & \textrm{max. exponent} \\
\hline\hline
generic-v10g40  & 10 &  40 &    52,131 & 29,987 \\
generic-v10g80  & 10 &  80 &   163,162 & 29,987 \\
generic-v10g120 & 10 & 120 &   411,997 & 29,991 \\
generic-v10g160 & 10 & 160 &   789,687 & 29,991 \\
generic-v10g200 & 10 & 200 & 1,245,139 & 29,991 \\
\hline
nongeneric-v10g100  & 10 &   100 &  19,442 & 10 \\
nongeneric-v10g150  & 10 &   150 &  52,781 & 10 \\
nongeneric-v10g200  & 10 &   200 &  79,003 & 10 \\
nongeneric-v10g400  & 10 &   400 & 193,638 & 10 \\
nongeneric-v10g600  & 10 &   600 & 318,716 & 10 \\
nongeneric-v10g800  & 10 &   800 & 435,881 & 10 \\
nongeneric-v10g1000 & 10 & 1,000 & 571,756 & 10 \\
\hline
squarefree-v20g100   & 20 &    100 &  3,990 & 1 \\
squarefree-v20g500   & 20 &    500 & 11,613 & 1 \\
squarefree-v20g2000  & 20 &  2,000 & 22,796 & 1 \\
squarefree-v20g4000  & 20 &  4,000 & 30,015 & 1 \\
squarefree-v20g6000  & 20 &  6,000 & 30,494 & 1 \\
squarefree-v20g8000  & 20 &  8,000 & 35,453 & 1 \\
squarefree-v20g10000 & 20 & 10,000 & 37,082 & 1 \\
\hline
J51        & 89  &   3,036 &          9 &  1 \\
J60        & 89  &   3,432 &         10 &  1 \\
smalldual  & 20  & 160,206 &         20 &  9 \\
frobn12d11 & 12  &  56,693 &  4,323,076 & 87 \\
frobn13d11 & 13  & 170,835 & 24,389,943 & 66 \\
k4         & 16  &      61 &        139 &  3 \\
k5         & 31  &  13,313 &     76,673 &  6 \\
model4vars & 16  &      20 &         64 &  2 \\
model5vars & 32  &     618 &      6,550 &  4 \\
tcyc5d25p  & 125 &   3,000 &     20,475 &  1 \\
tcyc5d30p  & 150 &   4,350 &     40,920 &  1 \\
\hline
\end{tabular}
\caption{Information about the test data.}
\label{fig:inputData}
\end{table}

\subsubsection*{Generation of random monomial ideals}
The random monomial ideals referred to below were generated using the
following algorithm, which depends on a parameter $N\in\N$. We start
out with the zero ideal. A random monomial is then generated by
pseudo-randomly generating each exponent within the range
$[0,N]$. Then this monomial is added as a minimal generator of the
ideal if it does not dominate or divide any of the previously added
minimal generators of the ideal. This process continues until the
ideal has the desired number of minimal generators. The random number
generator used was the standard C rand() function.

\subsubsection*{Description of the input data}

This list provides information on each test input.
\begin{description}

\item[generic] These ideals are nearly generic due to choosing
$N=30.000$.

\item[nongeneric] These ideals are non-generic due to choosing
$N=10$.

\item[square free] These ideals are square free due to choosing $N=1$.

\item[J51, J60]

These ideals were generated using the reverse engineering algorithm of
\cite{bionet}, and they were kindly provided by M. Paola Vera
Licona. They have the special features of having many variables, being
square free and having a small irreducible decomposition.

\item[smalldual]

This ideal has been generated as the Alexander dual of a random
monomial ideal with 20 minimal generators in 20 variables. Thus it has
many minimal generators and a small decomposition.

\item[t5d25p, t5d30p]

These ideals are from the computation of cyclic tropical polytopes,
and they have the special property of being generated by monomials of
the form $x_ix_j$ \cite{tropicalHull}. They were kindly provided by
Josephine Yu.

\item[k4, k5]

These ideals come with the program Monos \cite{monos} by R. Alexander
Milowski. They are involved in computing the integer programming gap
of a matrix \cite{igap}.

\item[model4vars, model5vars]

These ideals come from computations on algebraic statistical models,
and they were generated using the program 4ti2 \cite{4ti2} with the
help of Seth Sullivant.

\item[frobn12d11, frobn13d11] These ideals come from the computation
of the Frobenius number of respectively 12 and 13 random 11-digit
numbers \cite{bfrob}.
\end{description}

\subsection{Split Selection Strategies}
\label{ssec:strategyBench}

In this section we evaluate the split selection strategies described
in Section \ref{sec:selection}. Table \ref{fig:splitData} shows the
results.

The most immediate conclusion that can be drawn from Table
\ref{fig:splitData} is that label splits do well on ideals that are
somewhat generic, while they fare less well on square free ideals when
compared with pivot splits. It is a surprising contrast to this that
the MinLabel strategy is best able to deal with J60.

Table \ref{fig:splitData} also shows that the pivot strategies are
very similar on square free ideals. This is not surprising, as the
only valid pivots on such ideals have the form $x_i$, and the pivot
strategies all pick the same variable.

The final conclusion we will draw from Table \ref{fig:splitData} is
that the Median strategy is the best split selection strategy on these
ideals, so that is the strategy we will use in the rest of this
section. The Minimum strategy is a very close second.

\begin{table}
\centering
\begin{tabular}{|l||r|r|r|r|}
\hline
strategy & generic- & nongeneric- & squarefree- & \\
 & v10g200 & v10g400 & v20g10000 & J60\\
\hline\hline
MaxLabel & 13s & 13s  & 224s & 19s \\
MinLabel & 14s & 13s  & 203s &  2s \\
VarLabel & 18s & 13s  & 213s & 13s \\
\hline
Minimum  & 13s & 14s  & 19s  & 3s \\
\bf Median & \bf 12s & \bf 11s  & \bf 20s & \bf 3s \\
Maximum  & 35s & 43s  & 19s  & 3s \\
\hline
MinGen   & 59s & 201s & 19s  & 4s \\
Indep    & 13s & 12s  & 21s  & 3s \\
GCD      & 18s & 20s  & 19s  & 3s \\
\hline
\end{tabular}
\caption{Empirical comparison of split selection strategies.}
\label{fig:splitData}
\end{table}

\subsection{Empirical Comparison To Other Programs}
\label{ssec:otherBench}

In this section we compare our implementation in Frobby \cite{frobby}
of the Slice Algorithm to other programs that compute irreducible
decompositions. There are two well known fast algorithms for computing
irreducible decompositions of monomial ideals.

\begin{description}
\item[Alexander Dual \cite{alexdual, monoTree}] This algorithm uses
Alexander duality and intersection of ideals. Its advantage is speed
on highly non-generic ideals.

\item[Scarf Complex \cite{scarfComplex, monoTree}] This algorithm
enumerates the facets of the Scarf complex by walking from one facet
to adjacent ones. The advantage of the algorithm is speed for generic
ideals, while the drawback is that highly non-generic ideals lead to
high memory consumption and bad performance. This is because the
algorithm internally transforms the input ideal into a corresponding
generic ideal that can have a much larger decomposition.
\end{description}

We have benchmarked the following three programs.

\begin{description}
\item[Macaulay 2 version 1.0 \cite{M2}] Macaulay 2 incorporates an
  implementation of the Alexander Dual Algorithm. The time consuming
  parts of the algorithm are written in C++.

\item[Monos version 1.0 RC 2 \cite{monos}] This Java
  program\footnote{There are two different versions of Monos that have
    both been released as version 1.0. We are using the newest
    version, which is the version 1.0 RC2 that was released in 2007.}
  incorporates Alexander Milowski's implementation of both the
  Alexander Dual Algorithm and the Scarf Complex Algorithm.

\item[Frobby version 0.6 \cite{frobby}] This C++ program is our
  implementation of the Slice Algorithm.
\end{description}

How these programs compare depend on what kind of input is used, so we
use all the inputs described in Section \ref{ssec:dataBench} to get a
complete picture. In order to run these benchmarks in a reasonable
amount of time, we have allowed each program to run for one hour on
each input and no longer. Each program has been allowed to use 512 MB
of RAM and no more, not including the space used by other programs. We
use the abbreviation OOT for ``out of time'', OOM for ``out of
memory'' and RE for ``runtime error''.

The benchmarks have all been run on the same Linux machine with a 2.4
GHz Intel Celeron CPU. The reported time is the user time as measured
by the Unix command line utility ``time''.

All of the data can be seen on Table \ref{fig:benchData}. The data
shows that Frobby is faster than the other programs on all inputs
except for smalldual. This is because the Alexander Dual Algorithm
does very well on this kind of input, due to the decomposition being
very small compared to the number of minimal generators. The
decompositions of J51 and J60 are also small compared to the number of
minimal generators, though from the data not small enough to make the
Alexander Dual Algorithm win out.

It is clear from Table \ref{fig:benchData} that Macaulay 2 has the
fastest implementation of the Alexander Dual Algorithm when it does
not run out of memory. As expected, the Scarf Complex Algorithm beats
the Alexander Dual Algorithm on generic ideals, while the positions
are reversed on square free ideals.

As can be seen from Table \ref{fig:benchData}, the other programs
frequently run out of memory. In the case of Macaulay 2, this is
clearly in large part due to some implementation issue. However, the
issue of consuming large amounts of memory is fundamental to both the
Alexander Dual Algorithm and the Scarf Complex Algorithm, since it is
necessary for them to keep the entire decomposition in memory, and
these decompositions can be very large - see frobn13d11 as an
example. The Slice Algorithm does not have this issue.

An advantage of the Slice Algorithm is that the inner and outer slices
of a pivot split can be computed in parallel, making it simple to make
use of multiple processors. The Scarf Complex Algorithm is similarly
easy to parallelize, while the Alexander Dual Algorithm is not as
amenable to a parallel implementation. Although Frobby, Macaulay 2 and
Monos can make use of no more more than a single processor, multicore
systems are fast becoming ubiquitous. Algorithmic research and
implementations must adapt or risk wasting almost all of the available
processing power on a typical system. E.g. a non-parallel
implementation on an eight-way system will use only 13\% of the
available processing power.

\begin{table}
\centering
\begin{tabular}{|l||r|r|r|r|}
\hline
Input & Frobby & Macaulay2 & Monos & Monos \\
& & & (Alexander) & (Scarf) \\
\hline\hline
generic-v10g40 & $<$1s & 512s* & 1632s & 14s \\
generic-v10g80 & 1s & OOM & OOT & 82s \\
generic-v10g120 & 4s & OOM & OOT & 332s \\
generic-v10g160 & 8s & OOM & OOT & OOM \\
generic-v10g200 & 12s & OOM & OOT & OOM \\
\hline
nongeneric-v10g100 & $<$1s & 138s* & 770s & 191s \\
nongeneric-v10g150 & 1s & OOM & OOT & OOT \\
nongeneric-v10g200 & 1s & OOM & OOT & OOT \\
nongeneric-v10g400 & 4s & OOM & OOT & OOM \\
nongeneric-v10g600 & 8s & OOM & OOT & OOM \\
nongeneric-v10g800 & 11s & OOM & OOT & OOM \\
nongeneric-v10g1000 & 15s & OOM & OOT & OOM \\
\hline
squarefree-v20g100 & $<$1s & 17s & 27s & 1015s \\
squarefree-v20g500 & 1s & 80s & 608s & OOM \\
squarefree-v20g2000 & 4s & OOM & OOT & OOM \\
squarefree-v20g4000 & 9s & OOM & OOT & OOM \\
squarefree-v20g6000 & 13s & OOM & OOT & OOT \\
squarefree-v20g8000 & 19s & OOM & OOT & OOT \\
squarefree-v20g10000 & 21s & OOM & OOT & OOT \\
\hline
J51 & 2s & 8s & 6s & OOM \\
J60 & 3s & 10s & 7s & OOM \\
smalldual & 1961s & RE & 559s & RE \\
frobn12d11 & 285s & OOM & OOT & OOT \\
frobn13d11 & 2596s & RE & OOT & RE \\
k4 & $<$1s & 2s & 2s & 22s \\
k5 & 108s & OOM & OOT & OOM \\
model4vars & $<$1s & 1s & 1s & 2s \\
model5vars & 2s & OOM & 896s & OOM \\
tcyc5d25p & 7s & OOM & OOM & OOM \\
tcyc5d30p & 16s & OOM & OOT & OOM \\
\hline
\end{tabular}
\\ *: This time has been included despite using more than 512 MB of memory.
\caption{Empirical comparison of programs for irreducible
decomposition.}
\label{fig:benchData}
\end{table}

\subsection{The Bound Technique}
\label{ssec:boundBench}

In this section we examine the impact of using the bound technique
from Section \ref{sec:apps} to compute Frobenius numbers.

Table \ref{fig:boundData} displays the time needed to solve a
Frobenius problem IDP with and without using the bound technique for
some split selection strategies. We have included a new selection
strategy Frob that works as Median, except that it selects the
variable that maximizes the increase of the lower bound value on the
inner slice.

It is clear from Table \ref{fig:boundData} that the Frob and Median
split selection strategies are much better than the others for
computing Frobenius numbers, and that Frob is a bit better than
Median. We also see that applying the bound technique to the best
split selection strategy improves performance by a factor of between
two and three.

\begin{table}
\centering
\begin{tabular}{|l||r|r||r|r|}
\hline
strategy & frob-n11d11 & frob-n11d11 & frob-n12d11 & frob-n12d11 \\
 & without bound & using bound & without bound & using bound \\
\hline\hline
\bf Frob & \bf 66s & \bf 22s & \bf 204s & \bf 93s \\
Median & 76s & 35s & 256s & 147s \\
Maximum & 226s & 189s & 805s & 712s \\
Minimum & 731s & 761s & 3205s & 3388s \\
\hline
\end{tabular}
\caption{Empirical evaluation of the bound technique.}
\label{fig:boundData}
\end{table}

\bibliographystyle{utcaps}
\bibliography{references}

\end{document}